\documentclass{amsproc}

\newtheorem{Theorem}{Theorem}

\newtheorem{theorem}{Theorem}[section]
\newtheorem{lemma}[theorem]{Lemma}
\newtheorem{proposition}[theorem]{Proposition}
\newtheorem{corollary}[theorem]{Corollary}
\theoremstyle{definition}
\newtheorem{definition}[theorem]{Definition}

\theoremstyle{remark}

\numberwithin{equation}{section}

\newcommand{\abs}[1]{\lvert#1\rvert}
\newcommand{\Abs}[1]{\Vert\lvert#1\rvert\Vert}

\newcommand{\R}{{\mathbb R}}
\newcommand{\Z}{{\mathbb Z}}

\newcommand{\Int}{{{\rm Int}\,}}

\newcommand{\N}{{\mathbb N}}

\newcommand{\C}{{\mathbb C}}
\newcommand{\Fix}{{\rm Fix}}
\newcommand{\Moeb}{\mbox{M\"ob}_+(S^1_\C)}

\begin{document}

\title[A generic property of circle diffeomorphisms]
{A generic dimensional property
of the invariant measures for circle diffeomorphisms}


\author{Shigenori Matsumoto}
\address{Department of Mathematics, College of
Science and Technology, Nihon University, 1-8-14 Kanda, Surugadai,
Chiyoda-ku, Tokyo, 101-8308 Japan
}
\email{matsumo@math.cst.nihon-u.ac.jp
}
\thanks{The author is partially supported by Grant-in-Aid for
Scientific Research (C) No.\ 20540096.}


\date{\today }

\begin{abstract}
Given any Liouville number $\alpha$, it is shown that the nullity of
the Hausdorff dimension of the invariant measure is generic in
the space of the orientation preserving $C^\infty$ diffeomorphisms of
the circle with rotation number $\alpha$.
\end{abstract}

\maketitle

\section{Introduction}
Denote by $F$ the group of the orientation preserving $C^\infty$
diffeomorphisms of the cirlce. For $\alpha\in\R/\Z$, denote by $F_\alpha$
the subspace of $F$ consisting of all the diffeomorphisms whose rotation
numbers are $\alpha$, and by $O_\alpha$ the subspace of $F_\alpha$ 
of all the diffeomorphisms that are
 $C^\infty$ conjugate to $R_\alpha$, the rotation
by $\alpha$. 

In \cite{Y1}, J.-C. Yoccoz showed that $O_\alpha=F_\alpha$ if $\alpha$
is a non-Liouville number. Before than that, 
M. R. Herman (\cite{H}, Chapt.\ XI) 
had obtained the converse by showing
that for any Liouville number $\alpha$ the subspace $O_\alpha$ is meager in
$F_\alpha$.

For $f\in F_\alpha$, $\alpha$ irrational, denote by $\mu_f$ the unique
probability measure on $S^1$ which is invariant by $f$.
The properties of $\mu_f$ reflect the regularity of the conjugacy of
$f$ to $R_\alpha$. In \cite{S}, Victoria Sadovskaya improved 
the above result of M. R. Herman as follows.
For $d\in[0,1]$ define
$$S^d_\alpha=\{f\in F_\alpha\mid \dim_H(\mu_f)=d\},$$
where $\dim_H(\cdot)$ denotes the Hausdorff dimension of a measure.
She showed that for any Liouville number $\alpha$ and any $d\in[0,1]$, 
the set $S^d_\alpha$ is nonempty.
Notice that the Hausdorff dimension is an invariant of the equivalence
classes of measures, and therefore
$\dim_H(\mu_f)<1$ implies that $\mu_f$ is singular w.\ r.\ t.\ the
Lebesgue measure.

In \cite{Y2}, J.-C. Yoccoz showed the following theorem

\begin{theorem} \label{ty}
For any
irrational number the space
$O_\alpha$ is dense in $F_\alpha$ in the $C^\infty$ topology.
\end{theorem}

The proof of V. Sadovskaya is based on the method of fast approximation by
conjugacy with estimate, developed in \cite{FS}, and if it is slightly
modified it can be combined with the above theorem
to show that for any Liouville number $\alpha$ and
for any $d\in[0,1]$ the set $S^d_\alpha$ is $C^\infty$ dense in $F_\alpha$.

On the other hand M. R. Herman (\cite{H}, Prop I.8, p.\ 167)
showed that the set $S$ of the diffeomorphsim $f\in F$ such that $\mu_f$ is
singular {\em contains} a $G_\delta$ set in the $C^1$ topology in $F$.

These two results joined together {\em does not} imply immediately
that $S\cap F_\alpha$ is a dense $G_\delta$ set in the $C^r$ topology,
as pointed out to the author by Mostapha Benhenda.
The purpose of this paper is to settle down the situation.
In fact we get a bit more.

\begin{Theorem} \label{t}
For any Liouville number $\alpha$, the set $S_\alpha^0$
contains a countable intersection
of $C^0$ open and $C^\infty$ dense subsets of $F_\alpha$.
\end{Theorem}

\section{Preliminaries}
2.1. \ \ An irrational number $\alpha$ is called a {\em Liouville
     number}
if for any $N\in\N$ there is $p/q$ ($(p,q)=1$) such that 
$\abs{\alpha-p/q}<1/q^N$. We call $\alpha$ a {\em
lower Liouville number} if the above $p/q$ satisfies in addition
that  $p/q<\alpha$.\footnote{There are Liouville numbers which
look like non-Liouville, e.\ g.\ badly approximable, from one side.}
For any lower Louville number $\alpha$, $N\in\N$ and $\delta>0$ 
there is $p/q$ such
that $\abs{\alpha-p/q}<\delta/q^N$ and $p/q<\alpha$.

\bigskip
2.2. \ \ For a metric space $Z$ and $d>0$, the 
$d$-dimensional Hausdorff measure
$\nu^d(Z)$ is defined by
$$
\nu^d(Z)=\lim_{\varepsilon\to 0}\ \inf\{\sum^\infty_{i=1}
r_i^d\mid \cup_iB(x_i,r_i)=Z,\ r_i\leq\varepsilon\},
$$
where $B(x,r)$ denotes the open metric ball centered at $x$ of radius $r$.
The {\em Hausdorff dimension} of $Z$, denoted by $\dim_H(Z)$, is defined by
$$
\dim_H(Z)=\inf\{d\mid \nu^d(Z)=0\}=\sup\{d\mid\nu^d(Z)=\infty\}.
$$
 
The {\em lower box dimension} of $Z$, denoted by $\underline\dim_B(Z)$, 
is defined by
$$
\underline\dim_B(Z)=\underline\lim_{\varepsilon\to 0}\frac{\log N(\varepsilon, Z)}
{\log(1/\varepsilon)},
$$
where $N(\varepsilon,Z)$ denotes the minimal cardinality of
$\varepsilon$-dense
subsets of $Z$.

Let $X$ be a compact metric space, and  $\mu$ a probability measure
on $X$.
The Hausdorff dimension $\dim_H(\mu)$ and the lower box dimension 
$\underline\dim_B(\mu)$ of $\mu$ are defined respectively by

\begin{eqnarray*}
\dim_H(\mu)&=& \inf\{\dim_HZ\mid \mbox{$Z\subset X$ is measurable,}\ \mu(Z)=1\},\\
\underline\dim_B(\mu)&=&\lim_{\varepsilon\to 0}\ \inf\{\underline\dim_B(Z)
\mid \mbox{$Z\subset X$ is measurable,}\ \mu(Z)>1-\varepsilon\}.
\end{eqnarray*}

It is well known that
$$
\dim_H(\mu)\leq\underline\dim_B(\mu).
$$

\medskip
2.3. \ \ The proof of Theorem \ref{t} is by the method of fast approximation
by conjugacy with estimate. Let us prepare inequalities about
the derivatives of circle diffeomorphisms which are 
necessary for the estimate.

For a $C^\infty$ function $\varphi$ on $S^1$, we define 
as usual the $C^r$ norm
$\Vert\varphi\Vert_r$ ($0\leq r<\infty$) by
$$
\Vert\varphi\Vert_r=\max_{0\leq i\leq r}\sup_{x\in S^1}\abs{\varphi^{(i)}(x)}.$$

For $f,g\in F$, define
\begin{eqnarray*}
\Abs{f}_r&=&\max\{\Vert f-{\rm id}\Vert_r,\ \Vert f^{-1}-{\rm id}\Vert_r,1\},\\
d_r(f,g)&=&\max\{\Vert f-g\Vert_r,\ \Vert f^{-1}-g^{-1}\Vert_r\}.
\end{eqnarray*}

Since we include $1$ in the definition of $\Abs{f}_r$, we get the following
inequality from the Fa\`a di Bruno formula (\cite{H}, p.42 or \cite{S}).

\begin{lemma} \label{l1}
For $f,g\in F$ we have
$$
\Abs{fg}_r\leq C_1(r)\,\Abs{f}_r^r\,\Abs{g}_r^r,$$
where $C_1(r)$ is a positive constant depending only on $r$.
\qed
\end{lemma}

The following inequality can be found as Lemma 5.6 of \cite{FS} or
as Lemma 3.2 of \cite{S}.

\begin{lemma} \label{l2}
For $H\in F$ and $\alpha,\beta\in\R/\Z$,
$$
d_r(HR_\alpha H^{-1},HR_\beta H^{-1})\leq
C_2(r)\,\Abs{H}_{r+1}^{r+1}\,\abs{\alpha-\beta},
$$
where $C_2(r)$ is a positive constant depending only on $r$.
\qed
\end{lemma}

For $q\in\N$, denote by $\pi_q:S^1\to S^1$ the $q$-fold covering map.
Simple computation shows:

\begin{lemma} \label{l3}
Let $h$ be a lift of $k\in F$ by $\pi_q$ and assume $\Fix(h)\neq\emptyset$. 
Then we have
$$
\Abs{h}_r\leq\Abs{k}_r\,q^{r-1}.$$
\qed
\end{lemma}

\smallskip
2.4. \ \ Here we prepare necessary facts about Moebius transformations
on the circle. Let
$$
S^1_\C=\{z\in\C\mid \abs{z}=1\},$$
and $\Moeb$  the group of the orientation preserving
Moebius transformations of $\C$ which leaves $S^1_\C$ invariant.
We identify $S^1_\C$ with the circle $S^1=\R/\Z$ in a standard way.
For $k\in\Moeb$, the diffeomorphism of $S^1$ corresponding
to $k$ is denoted by $\hat k$. Define the {\em expanding interval}
$\mathcal I(\hat k)$ of $\hat k$ by
$$
\mathcal I(\hat k)=\{x\in S^1\mid \hat k'(x)\geq 1\}.
$$
Then the inverse formula of the derivatives shows that
$$
\hat k(\mathcal I(\hat k))=S^1\setminus\Int\mathcal I(\hat k^{-1}).$$
Denote by $\rho(\hat k)$ the radius of $\mathcal I(\hat k)$.
Notice that $\rho(\hat k)=\rho(\hat k^{-1})$. 

For $1/2\leq a<1$, define $k_a\in\Moeb$ by
$$
k_a(z)=\frac{z+a}{az+1}.$$
The transformation $k_a$ is hyperbolic with an attractor $z=1$
and a repellor $z=-1$. Notice that $\abs{k_a'(-1)}\nearrow\infty$
as $a\nearrow 1$. The corresponding diffeomorphism $\hat k_a$
has an attractor at $x=0$ and a repellor at $x=1/2$, and
$\rho(\hat k_a)\searrow 0$ as $a\nearrow 1$.

\begin{lemma} \label{l4}
There is a constant $C_3(r)>0$ depending only on $r$ such that
for any $1/2\leq a<1$,
$$
\Abs{\hat k_a}_r\leq C_3(r)\rho(\hat k_a)^{-2r}.$$
\end{lemma}

{\sc Proof}. First of all $\rho(\hat k_a)$ is proportional to
the radius of the isometric circle of $k_a$,
$\{z\in\C\mid\abs{k_a'(z)}=1\}$,
and the latter can easily be computed using the expression
\begin{equation} \label{e1}
k_a'(z)=\frac{1-a^2}{(az+1)^2}.
\end{equation}
It follows that there is a constant $c>0$
such that
$$\rho(\hat k_a)\leq c(1-a)^{1/2},\ \ \ 1/2\leq a<1.$$

For $k_a$, looked upon as a map from $S^1_\C$ to $S^1_\C$, the real
$r$-th derivative w.\ r.\ t.\ the angle coordinate 
is denoted by $D^rk_a$, while 
$\varphi'$ denotes the complex derivative of a holomorphic map $\varphi$.
It suffices to show for any $r$ and $z\in S^1_\C$,
\begin{equation} \label{e2}
\abs{D^rk_a(z)}\leq c_3'(r)(1-a)^{-r}.
\end{equation}

For $r=1$, this follows immediately from (\ref{e1}) since $Dk_a=\abs{k_a'}$.

Now $Dk_a$ extends to a holomorphic function on a neighbourhood
of
$S^1_\C$.
Since $\arg(k_a'(z))=\arg(k_a(z)/z)$ and $\abs{k_a(z)/z}=1$
for $z\in S^1_\C$,
we have 
\begin{equation*} \label{e3}
Dk_a(z)=k_a'(z)z/k_a(z)=(\frac{1}{z+a}-\frac{a}{az+1})z.
\end{equation*}

It follows that
$$
D^2k_a=\abs{(Dk_a)'},
$$
where
$$
(Dk_a)'(z)=\frac{P_1}{(z+a)^2}+\frac{Q_1}{(az+1)^2},$$
and $P_1$ and $Q_1$ are polynomials in $z$ and $a$, showing
(\ref{e2}) for $r=2$. 

Now since $Dk_a$ is real valued on $S^1_\C$,
its derivative along the direction tangent to $S^1_\C$ is real.
Therefore $D^2k_a$ extends to a holomorphic function as
$$D^2k_a(z)=(Dk_a)'(z)iz.$$
This shows that

$$
D^3k_a=\abs{(D^2k_a)'},$$
where
$$
(D^2k_a)'(z)=\frac{P_2}{(z+a)^3}+\frac{Q_2}{(az+1)^3},
$$
showing (\ref{e2}) for $r=3$. 

The last argument for $r=3$ can be applied
for any $r\geq 4$, completing the proof of the lemma.
\qed

\section{The $G_\delta$ set}
In the rest of the paper we choose an arbitrary
Liouville number $\alpha$ and fix it once and for all.
Let us assume that $\alpha$ is a lower
Liouville number (See 2.1.), the other case being dealt with similarly.
In this section we define a $G_\delta$ set $B$ of
$F_\alpha$ in the $C^0$ topology, and show that any $f\in B$ satisfies
that $\dim_H(\mu_f)=0$. Notice that by the lower Liouville property,
for $p/q$ well approximating $\alpha$, the iterate $R_\alpha^q$ has
rotation number $q\alpha-p$, a very small positive number.

\begin{definition} \label{d1}
For any $n\in\N$, we define $B_n \subset F_\alpha$ to be
the subset consisting of those $f$ which satisfy the following
condition.

There exist integers $q_n=q_n(f)>n$, $l_n=l_n(f)>0$ and points
$c_i=c_i^n(f)$, $d_i=d_i^n(f)$ of $S^1$ ($0\leq i\leq q_n$) with
$c_{q_n}=c_0$ and $d_{q_n}=q_0$,
satisfying the following properties.

\smallskip
(\ref{d1}.a) $c_1<d_1<c_2<d_2<\cdots<c_{q_n}<d_{q_n}<c_1$
in the cyclic order.

\smallskip
(\ref{d1}.b) $\max_i(d_i-c_i)<q_n^{-n}$.

\smallskip
(\ref{d1}.c) $\max_i(d_i-c_i)<n^{-1}\min_i(c_{i+1}-d_i)$.

\smallskip
(\ref{d1}.d) $f^{kq_n}(c_i)\in(c_i,d_i)$ for any $1\leq k\leq 2^nl_n$.

\smallskip
(\ref{d1}.e) $f^{l_n q_n}(d_i)\not\in[d_i,c_{i+1}]$.
\end{definition}

Clearly $B_n$ is $C^0$ open in $F_\alpha$, and therefore
their intersection $B=\cap_n B_n$ is a $G^\delta$-set.

The following lemma follows from the flexibility of Definition
\ref{d1}, e.\ g.\ (\ref{d1}.c).

\begin{lemma} \label{l31}
For $h\in F$, we have $hBh^{-1}=B$.
\qed
\end{lemma}

\begin{lemma} \label{l32}
If $f\in B$, then $\dim_H(\mu_f)=0$.
\end{lemma}

{\sc Proof}. It suffices to show that $\underline\dim_B(\mu_f)=0$.
(See Paragraph 2.2.) \ Choose an arbitrary $f\in B$. Below we depress the
notations
as
$$q_n=q_n(f),\ l_n=l_n(f),\   c_i=c_i^n(f),\  d_i=d_i^n(f).$$
Define
$$I^{(n)}=\bigcup_{i=1}^{q_n}[c_i,d_i].$$
By (\ref{d1}.d) and (\ref{d1}.e), we have
$$ \#\{i\mid f^{iq_n}(x)\in I^{(n)},\  1\leq i\leq 2^nl_n\}\geq (2^n-1)l_n,
\ \ \forall x\in S^1.
$$
This implies $$\mu_f(I^{(n)})\geq1-2^{-n}.$$
Define a closed set $C_m=\cap_{n\geq m}I^{(n)}.$ Then since
$$
\mu_f(C_m)\geq 1-2^{-m+1},$$
it suffices to show that $\underline\dim_B(C_m)=0$.
Choose $\varepsilon=\max_i(d_i-c_i)$. Then by (\ref{d1}.c) we have for 
$m\geq2$,
$$N(\varepsilon, C_m)\leq N(\varepsilon, I^{(n)})=q_n, \ \ \forall n_\geq m.$$
We also have $\varepsilon\leq q_n^{-n}$ by (\ref{d1}.b), showing that
$$
\frac{\log N(\varepsilon, C_m)}{\log(1/\varepsilon)}\leq n^{-1}.$$
Since $n$ is an arbitrary integer $\geq m$, we have shown
that $\underline\dim_B(C_m)=0$, as is required.
\qed

\section{proof}
The purpose of this section is to show:

\begin{proposition} \label{p1}
For any $r\in\N$, there is $f\in B$ such that 
$d_r(f,R_\alpha)<2^{-r}$.
\end{proposition}

The proposition asserts that $R_\alpha$ belongs to the $C^\infty$ closure
of $B$, which implies that $B$ is $C^\infty$ dense in $F_\alpha$,
by virtue of Theorem \ref{ty} and Lemma \ref{l31}.
This, together with the fact that
$B$ is a $G_\delta$ set of $F_\alpha$
in the $C^0$ topology, completes the proof of Theorem \ref{t}.

Our overall strategy of the proof of Proposition \ref{p1} is as follows.
Since $\alpha$ is lower Liouville we can choose a sequence of
rationals $\alpha_n=p_n/q_n$ well approximating $\alpha$
so that $\alpha_n\nearrow\alpha$.
We shall construct a diffeomorphism $h_n\in F$ 
commuting with $R_{\alpha_n}$, and set
$$
H_n=h_1h_2\cdots h_n,\ \ f_n=H_nR_{\alpha_{n+1}}H_n^{-1}.$$
We are going to show that $f_n$'s converge to $f\in B$.
The commutativity condition above is quite useful when we estimates
the norm of the functions
$$
f_{n-1}^{\pm 1}-f_{n}^{\pm 1}=H_nR_{\alpha_n}^{{\pm 1}}H_n^{-1}
-H_nR_{\alpha_{n+1}}^{{\pm 1}}H_n^{-1},$$
by virtue of Lemma \ref{l2}.

Now a concrete construction gets started. Choose $\alpha_1=p_1/q_1$
so that

\smallskip
(A) $0<\alpha-\alpha_1<2^{-(r+1)}$,

\smallskip
\noindent
and let $f_0=R_{\alpha_1}$.

Set $h_1$ to be the lift of $\hat k_{a_1}$ by the  
cyclic $q_1$-covering such that $\Fix(\hat k_{a_1})\neq\emptyset$,
where $a_1\in[1/2,1)$ satisfies

\smallskip
(${\rm B}_1$) $\rho(\hat k_{a_1})=q_1^{-2}$.

\smallskip
\noindent
See Paragraph 2.4 for these definitions.
Notice that $h_1R_{\alpha_1}h_1^{-1}=R_{\alpha_1}$.

Assume we already defined $\alpha_i$ and $h_i$ for $1\leq i\leq n-1$.
Define
$$
L_{n-1}=\max\{\Vert H_{n-1}'\Vert_0, \Vert(H_{n-1}^{-1})'\Vert_0\}.$$
Choose $\alpha_n=p_n/q_n$ which satisfies (C), (D) and (E) below.
Such $\alpha_n$ exits since $\alpha$ is a lower Liouville number.
The constants $C_i(\cdot)$, $i=1,2,3$, are from the lemmata in Sect.\ 2.

\smallskip
(C) $0<\alpha-\alpha_n<\delta/q_n^N$, where

$$
\delta=2^{-(n+r+1)}C_2(n+r)^{-1}C_1(n+r+1)^{-(n+r+1)}
C_3(n+r+1)^{-(n+r+1)^2}\Abs{H_{n-1}}_{n+r+1}^{-(n+r+1)^2},$$
$$N=(2n+3)(n+r+1)^3.$$

\smallskip
(D) $q_n>2^nL_{n-1}.$

\smallskip
(E) $q_n>2^{n+5}q_{n-1}.$

\smallskip
Finally set $h_n$ to be the lift of $\hat k_{a_n}$ by the $q_n$ covering
with $\Fix(h_n)\neq\emptyset,$, where $1/2\leq a_n<1$ is chosen 
such that

\smallskip
(B) $\rho(\hat k_{a_n})=q_n^{-(n+1)}$.

\smallskip
Notice that $h_nR_{\alpha_n}h_n^{-1}=R_{\alpha_n}$
and that $\Vert h_n^{\pm1}-{\rm id}\Vert_0\leq 2^{-1}q_n^{-1}$.
\begin{lemma} \label{l41}
We have $d_{n+r}(f_{n-1},f_n)<2^{-(n+r+1)}$ for $n\geq1$.
\end{lemma}

{\sc Proof}. The proof is a routine calculation using the lemmata
 in Sect.\ 2 and condition (C). 
Just notice that $f_{n-1}=H_nR_{\alpha_n}H_n^{-1}$,
while $f_{n}=H_nR_{\alpha_{n+1}}H_n^{-1}$,
and that $0<\alpha_{n+1}-\alpha_n<\alpha-\alpha_n$.
\qed

\begin{corollary} \label{c1}
The limit $f=\lim_{n\to\infty}f_n$ is a $C^\infty$ diffeomorphism and
$d_r(f,R_\alpha)\leq2^{-r}$.
\end{corollary}

{\sc Proof}. The latter assertion is obtained from (A) and the following
estimate.
$$
d_r(f,R_{\alpha_1})\leq\sum_{n=1}^\infty d_r(f_{n-1},f_n)
\leq \sum_{n=1}^\infty 2^{-(n+r+1)}\leq 2^{-(r+1)}.
$$
\qed

\begin{lemma} \label{l42}
There exists a homeomorphism $H$ of $S^1$ such that 
$d_0(H_n,H)\to 0$.
\end{lemma}

{\sc Proof}. First we have by (E)
$$\Vert H_n^{-1}-H_{n-1}^{-1}\Vert_0\leq\Vert h_n^{-1}-{\rm id}\Vert_0\leq
2^{-1}q_n^{-1}
\leq 2^{-n}.$$
On the other hand by (D)
$$
\Vert H_n-H_{n-1}\Vert_0\leq L_{n-1}\Vert h_n-{\rm id}\Vert_0\leq
L_{n-1}2^{-1}q_n^{-1}\leq 2^{-(n+1)},
$$
showing the lemma. \qed

\medskip
It follows from Lemma \ref{l42} that $f=HR_{\alpha}H^{-1}$ and in particular $f\in
F_\alpha$.

\bigskip
{\em In what follows, we fix $n\in\N$ once and for all and will show that
$f\in B_n$}.
First of all let us study the dynamics of $h_n$ in details.
Recall that $h_n$ is a lift of $\hat k_{a_n}$ by the $q_n$-fold
covering. So $h_n$ has $q_n$ repelling fixed points
and $q_n$ attracting fixed points.

The expanding interval $\mathcal I(\hat k_n)$ of $\hat k_n$ (See 2.4.)
is centered at $1/2$ and has length $2q_n^{-(n+1)}$ by (B), and
$\mathcal I(\hat k_n^{-1})$ is the interval centered at $0$ of the
same length. Recall the dynamics of $\hat k_{a_n}^{-1}$:
$$
\hat k_{a_n}^{-1}(\mathcal I(\hat k_{a_n}^{-1}))
= S^1\setminus \Int \mathcal I(\hat k_{a_n}).
$$
Let $[c_i',d_i']$, $1\leq i\leq q_n$ be the lift of 
$\mathcal I(\hat k_{a_n}^{-1})$, located in this order in $S^1$.
Their lengths $d'_i-c'_i$ is very small compared with
$c'_{i+1}-c'_i=q_n^{-1}$. In fact by (B)
$$
d'_i-c'_i=2q_n^{-(n+2)}.
$$

The intervals $[h_n^{-1}d'_{i}, h_n^{-1}c'_{i+1}]$ are lifts
of $\mathcal I(\hat k_{a_n})$, and has the same length $2q_n^{-(n+2)}$.
Since by (E)
$$2q_n^{-(n+2)}<2^{-(n+3)}q_n^{-1},$$ 
 we have
\begin{equation} \label{e41}
0<h_n^{-1}c'_{i+1}-h_n^{-1}d'_{i}<2^{-(n+3)}q_n^{-1}.
\end{equation}
Put
$$
H_{(n)}=\lim_{k\to\infty}h_nh_{n+1}\cdots h_{n+k},\ \ \
f_{(n)}=H_{(n)}R_\alpha H_{(n)}^{-1},$$
$$
c''_i=H_{(n)}^{-1}c'_i,\ \ \ d''_i=H_{(n)}^{-1}d'_i.$$
Then we have
\begin{eqnarray} 
\abs{c''_i-h_n^{-1}c'_i} & \leq & 2^{-(n+5)}q_n^{-1}, \label{ee}
\\
\abs{d''_i-h_n^{-1}d'_i} & \leq & 2^{-(n+5)}q_n^{-1}. \label{eee}
\end{eqnarray}

In fact, (\ref{ee}) can be shown by
$$\abs{c''_i-h_n^{-1}c'_i}\leq\Vert H_{(n+1)}^{-1}-{\rm id}\Vert_0
\leq \sum_{i=1}^\infty\Vert h_{n+i}^{-1}-{\rm id}\Vert_0
$$
$$
\leq \sum_{i=1}^\infty 2^{-1}q_{n+i}^{-1}\leq
\sum_{i=1}^\infty 2^{-(n+i+5)}q_n^{-1}=2^{-(n+5)}q_n^{-1},
$$
where we have used (E) in the last inequality.

From (\ref{e41}), (\ref{ee}) and (\ref{eee}), we obtain that
\begin{equation} \label{eeee}
0<c''_{i+1}-d''_i< 2^{-(n+2)}q_n^{-1}.
\end{equation}
On the other hand we have
\begin{equation} \label{eeeee}
d''_i-c''_i>2^{-1}q_n^{-1}.
\end{equation}
In fact
$$d''_i-c''_i=(c''_{i+1}-c''_{i})-(c''_{i+1}-d''_i)
>q_n^{-1}-2\cdot 2^{-(n+5)}q_n^{-1}-2^{-(n+2)}q_n^{-1}>2^{-1}q_n^{-1}.
$$

Now the rotation number of $R_\alpha^{q_n}$ is $q_n\alpha-p_n$,
a very small positive number. Let us estimate how long the orbit 
by $R_\alpha^{q_n}$ of $c''_i$ stays in the interval $(c''_i,d''_i)$.
Let $m_n$ be the largest integer such that 
$$R_\alpha^{kq_n}c''_i\in(c''_i,d''_i), \ \ \mbox{if} 
\ \ 1\leq k\leq m_n.$$
Then we have by (\ref{eeeee})
$$
m_n=\lfloor (d''_i-c''_i)(q_n\alpha-p_n)^{-1}\rfloor
\geq 2^{-1}q_n^{-1}(q_n\alpha-p_n)^{-1}-1.$$
Next estimate
how quickly the orbit of $d''_i$ exits $[d''_i,c''_{i+1}]$.
Let $l_n$ be the smallest positive integer such that 
$$R_\alpha^{l_nq_n}d''_i\not\in[d''_i,c''_{i+1}].$$
 Then it follows
from (\ref{eeee}) that
$$
l_n=\lfloor (c''_{i+1}-d''_i)(q_n\alpha-p_n)^{-1}\rfloor
+1\leq 2^{-(n+2)}q_n^{-1}(q_n\alpha-p_n)^{-1}+1.$$
By (C) the number $q_n^{-1}(q_n\alpha-p_n)^{-1}$ is
sufficiently big, and we have
\begin{equation} \label{e42}
m_n\geq2^nl_n.
\end{equation}

Now consider $f=HR_\alpha H^{-1}.$
Let $c_i=Hc''_i$ and $d_i=Hd''_i$.
Then $m_n$ is
the largest integer such that $f^{kq_n}c_i\in(c_i,d_i)$
if $1\leq k\leq m_n$ and $l_n$ the smallest positive integer such that 
$f^{l_nq_n}d_i\not\in[d_i,c_{i+1}]$. Thus (\ref{e42}) implies
(\ref{d1}.d) and (\ref{d1}.e).

Finally recall that
$$
f=HR_\alpha H^{-1}=H_{n-1}f_{(n)}H_{n-1}^{-1}, \ \ \
c_i=H_{n-1}c'_i, \ \ \ d_i=H_{n-1}d'_i.$$ 
As for (\ref{d1}.b) we have
\begin{equation} \label{e43}
d_i-c_i\leq L_{n-1}(d'_i-c'_i)=L_{n-1}2q_n^{-(n+2)}.
\end{equation}

Now (D) implies that 
$$L_{n-1}2q_n^{-(n+2)}\leq q_n^{-n}.$$
This shows (\ref{d1}.b).

For (\ref{d1}.c), we have 
\begin{equation} \label{e44}
c_{i+1}-d_i\geq L_{n-1}^{-1}(c'_{i+1}-d'_{i})\geq L_{n-1}^{-1}
2^{-1}q_n^{-1}.
\end{equation}

Also (D) implies 
\begin{equation} \label{e45}
4nL_{n-1}^2\leq q_n^{n+1}.
\end{equation}
Simple computation shows that (\ref{e43}), (\ref{e44}) and (\ref{e45})
implies the condition (\ref{d1}.c).
Now we are done with the proof that $f\in B_n$.
Since $n$ is arbitrary, this shows that $f\in B$, completing
the proof of Proposition \ref{p1}.

\end{document}